\documentclass[11pt]{article}
\usepackage{amsmath, amscd, amsthm} 
\usepackage{amssymb, amsfonts}
\input xy \xyoption{all}
\CompileMatrices

\def\proof{\par\medskip\noindent {\bf Proof: }}
\def\qed{\hfill $\Box$ \medskip \par}
\def\Box{\framebox[10pt]{\rule{0pt}{3pt}}}

\ifx\ThmNum\undefined 
   \newtheorem{theorem}{Theorem}[section]
   
\else 
   \newtheorem{theorem}{Theorem} 
   
\fi

\newcommand{\Z}{\mathbb{Z}}
\newcommand{\R}{\mathbb{R}}
\newcommand{\C}{\mathbb{C}}
\newcommand{\Q}{\mathbb{Q}}
\newcommand{\F}{\mathbb{F}}
\newcommand{\SH}{\mathcal{S}\mathcal{H}}

\newcommand{\KMW}{K^{MW}}
\newcommand{\KM}{K^{M}}
\newcommand{\kM}{k^{M}}

\newtheorem{proposition}[theorem]{Proposition}
\newtheorem{lemma}[theorem]{Lemma}
\newtheorem{definition}[theorem]{Definition}
\newtheorem{corollary}[theorem]{Corollary}

\newtheorem{question}[theorem]{Question}

\newtheorem{remark}[theorem]{Remark}

\title{Some comments on motivic nilpotence}
\author{Jens Hornbostel}
\date{January 20, 2017}

\begin{document}

\maketitle

\begin{abstract}
We discuss some results and conjectures related
to the existence of the non-nilpotent motivic maps
$\eta$ and $\mu_9$. To this purpose, we establish
a theory of power operations for motivic
$H_{\infty}$-spectra. Using this, we show that
the naive motivic analogue of the unstable
Kahn-Priddy theorem fails.  
Over the complex numbers, we show
that the motivic $T$-spectrum
$S[\eta^{-1},\mu_9^{-1}]$
is closely related to higher Witt groups,
where $S$ is the motivic sphere spectrum and
$\eta$ and $\mu_9$ are explicit
elements in $\pi_{**}(S)$. 
\end{abstract}

\section{Introduction}

Let us fix a base field $k$, and consider the motivic stable
homotopy category $\SH(k)$ of Morel-Voevodsky.
In this article, we investigate possible motivic generalizations of  
three famous theorems in classical homotopy theory:
Nishida's nilpotence theorem (see \cite{N} and also \cite{M1},
\cite{M2}), the fact that complex cobordism $MU_*$ detects 
nilpotent self maps \cite{DHS}
and the classification of thick ideals in $\SH_{fin}$
resp. its $p$-localizations via Morava-$K$-theories 
$K(n)$ \cite{HS}.
These theorems lay the foundations of the chromatic
approach to stable homotopy theory. 
The theorems are closely related classically,
and motivic counterparts would be as well:
We have a decomposition of Bousfield classes
$<MU_{(p)}>= \bigvee_{n \geq 0} <K(n)>$, and in the motivic
case Joachimi \cite[Theorem 9.5.1]{J} has shown a slightly weaker decomposition
theorem over $\C$ for $p$ odd. Moreover,
any non-nilpotent map $f \in 
\pi_{**}(S)$ produces a thick ideal $thickid(Cf)$
which consists precisely of those objects $X$ for which
$f^n \wedge id_X$ is zero for $n>>0$, as 
shown in \cite[Theorem 2.15]{Ba}. Besides elements in bidegree
$(0,0)$, the first
non-nilpotent motivic self map is the 
{\em motivic Hopf map} $\eta \in \pi_{1,1}(S)$.
It has been studied by Morel who has shown
that unlike its topological
counterpart it is {\em not} nilpotent,
but is not detected by algebraic
cobordism ${\bf MGL}$ either \cite{Mo}.
Hence $\eta$ shows that at least two of three
above theorems do not hold motivically.
Applying the above result to $\eta$, 
one obtains a thick ideal $thickid(C\eta)$ which 
Joachimi \cite[chapter 7]{J} shows to be different
from the ones arising from classical or $\Z/2$-equivariant
stable homotopy theory, thus disproving a motivic 
version of the third theorem above. Recall 
that one may define motivic Morava $K$-theories
by killing elements in ${\bf MGL}_{**}$, as done first by
Borghesi and Hu. See also \cite{Ho2}
for chromatic versions of Johnson-Wilson spectra
and \cite{J} for more Bousfield class
decompositions over $\C$.

\medskip

Similar to the study of localizations of the topological sphere spectrum
with respect to given integers in $\pi_0(S)$,
one is naturally led to study the $\eta$-localization of the motivic sphere
spectrum $S$ over $k$. For $k=\C$, Andrews and Miller \cite{AM}
recently proved the following, confirming
a conjecture of \cite{GuI}:
\begin{theorem}(Andrews-Miller)\label{AMbeautiful}
For $k=\C$, the motivic
ANSS induces an isomorphism
$$\pi_{**}(S)[\eta^{-1}] \cong \F_2[\eta^{\pm 1}, \sigma, \mu_9]/
(\eta \sigma^2),$$ with $\mu_9  \in \, <8 \sigma,2,\eta> \, \subset \pi_{9,5}(S)$
and $\sigma \in \pi_{7,4}(S)$ the motivic Hopf element of \cite{DI3}.
\end{theorem}

Our notation follows the usual motivic grading convention in which the simplicial one-dimensional sphere is 
\(S^{1,0}\) and the punctured affine line is \(S^{1,1}\).

\medskip

The element $\mu_9$ considered in \cite{AM}
is originally defined only after motivic $2$-completion.
However, it can be shown to exist integrally, see Lemma \ref{mu9int} below.
It is detected by $Ph_1$ in the motivic ASS and satisfies
$2 \mu_9$=0 \cite[p.1010]{DI2}, \cite[table 8]{Is}.
We refer to \cite[table 23]{Is} for the indeterminacy
of $<8 \sigma,2,\eta>$. Let us emphasize
the importance of the huge amount of computaions by Isaksen
(and coauthors) at the prime $2$, see in particular \cite{Is}. 

We now fix the element $\mu_9
\in \pi_{9,5}(S)$ which is detected by the motivic
$\alpha_5$ in the motivic ANSS, see \cite[section 7]{AM}. 
More than a decade after Morel's study
of $\eta$, the existence of $\mu_9$ yields the second example
of an element -- at least over $\C$ -- contradicting
the naive motivic analogue of the classical three
nilpotence theorems above:

\begin{corollary}
Let $k=\C$.

(i) The above element $\mu_9 \in <8 \sigma,2,\eta>$ is not nilpotent
in $\pi_{**}(S)$.

(ii) The motivic algebraic cobordism spectrum
${\bf MGL}$ does not detect $\mu_9$,
i.e. ${\bf MGL}(\mu_9)=0$.

(iii) The thick ideal $thickid(C\mu_9)$ in $\SH(\C)$
is ``new'', that is it is not induced by
some thick ideal of $\SH$.
\end{corollary}
\proof
(i) Follows immediatley from Theorem \ref{AMbeautiful}.
(ii) This is an easy consequence of the computation
in \cite[Theorem 7]{HKO1}. (As Andrews pointed out to me, this also follows
because $\mu_9$ has Novikov filtration one.)
(iii) This is similar to \cite[Proposition 7.1.4 (2),(3)]{J},
using \cite[Corollary 2.15]{Ba} (which implies that
$thickid(C\mu_9)$ is strictly smaller than $\SH(\C)_{fin}$)
and that the complex realization of $\mu_9$ is nilpotent.
\qed

The description of $\mu_9$ as a Toda bracket 
$<8 \sigma,2,\eta>$ of motivic Hopf 
elements together with the fact $(1 - \epsilon)\eta=0$
for all fields shows that it exists over other
fields as well, provided that
$2^t \cdot (1 - \epsilon)^{4-t} \cdot \sigma=0$
for some $t$. As Isaksen points out, computations in
\cite{DI4} and subsequent work imply that this is non-zero for
$k=\R$. Hence $\mu_9$ is not defined for subfields
$k \subset \R$. If $\mu_9$ is defined for some 
subfield $k \subset \C$, then it is
a non-zero non-nilpotent element over $k$. 
(It remains to compare $thickid(C\eta)$ and $thickid(C\mu_9)$,
of course.) 

After the appearance of $\mu_9$, more non-nilpotent self
maps have been recently discovered by Isaksen and his coauthors.
We refer the reader to the beginning of the next section
for a short discussion of this work in progress.
Let us however mention the following here: Andrews conjectures 
that $\eta=w_0$ is just the beginning of a new chromatic motivic
family $w_0,w_1,w_2...$ for $k=\C$ at the prime $2$.
He has constructed self maps $w_1^4:S/\eta \to S/\eta$
and $w_2^{32}$. Inspired by his construction, Gheorge
has constructed a new infinite family $K(w_n)$ of
motivic Morava-$K$-theory spectra.

The very recent preprint \cite{Th}
computes the homogenous spectrum of $K^{MW}_*(k)$.
This together with work in progress by Heller and Ormsby
(the preprint is now available, see \cite{HO}) 
showing that the graded map \cite{Ba} to the homogenous 
spectrum of $K_*^{MW}(k)$ is surjective leads
to a refinement of
the results above about thick (prime) ideals in $\SH(k)$
including $C\eta$, but not $C\mu_9$.
(The case of a finite base field has been considered
slightly earlier by \cite{K}.)

\medskip

This introductory section 
contains mostly recollections of known
recent results and some rather easy consequences
thereof. The main results appear in the following section 2.
There is a very close and precise relationship
between the failure of motivic Nishida nilpotence and 
the motivic unstable Kahn-Priddy theorem, see Theorem \ref{main}.
This has very explicit consequences, see Corollary \ref{mainappl}. 
For instance, for $k=\C$ the motivic Hopf element $\eta \in 
\pi_{9,3}(S^{8,2})$ does not lift to $\pi_{9,3}((E \Z/2)_+ 
\wedge_{\Z/2}S^{4,1})$, although the topological Hopf 
$\eta_{top}$ does lift to $\pi_{9}((E \Z/2)_+ \wedge_{\Z/2}S^{4})$. 
In other words, a certain homotopical symmetry
of $\eta_{top}$ does not lift to the motivic
Hopf map $\eta$.

The proofs of these results rely heavily on the study of 
{\em extended powers} and {\em power operations}
for {\em motivic} $H_{\infty}$-spectra, which generalizes the
classical work of J-P. May et al and represents a topic 
of interest for its own sake.

Finally, over the complex numbers there is a beautiful
relationship between the motivic sphere spectrum and the
spectrum ${\bf KT}$ representing Witt groups (see Theorem \ref{main2}):

\begin{theorem}
For $k=\C$, the unit map of the 
hermitian $K$-theory spectrum
induces an epimorphism 
$$\pi_{**}(S[\eta^{-1},\mu_9^{-1}]) \to 
\pi_{**}{\bf KT}$$
with the kernel being another copy of $\pi_{**}{\bf KT}$. 
\end{theorem}

This theorem is an incarnation of the following general
philosophy: quadratic forms and hermitian $K$-groups detect a 
lot of interesting motivic homotopy theory not visible by cohomology or
classical homotopy theory. (In the classical case,
this is reflected e.g. by the Moebius stripe detecting the 
topological Hopf map.) This philosophy has been confirmed 
before e.g. by \cite{Mo}, \cite{RSO2} and \cite{ALP}. Our
theorem may be considered as an integral refinement
of the latter over the complex numbers.

The appendix of Marcus Zibrowius provides a complete
answer to motivic Nishida nilpotence in simplicial
degree zero, that is for Milnor-Witt $K$-theory.
The main result is the following 
(see Proposition~\ref{app-prop}):
  
\begin{theorem}(Zibrowius)\label{app-prop}
Nishida nilpotence holds in all non-negative degrees of the Milnor-Witt 
K-ring \(\KMW_*(F)\).
Nishida nilpotence holds in all degrees of \(\KMW_*(F)\) if and only 
if \(F\) is formally real.
\end{theorem}

\medskip

We thank Tom Bachmann, Bert Guillou, Dan Isaksen, Karlheinz Knapp, Sean Tilson
and Marcus Zibrowius for useful comments, and Dan Isaksen also for informing us about 
some of his computations over $\R$.

\section{On the failure of motivic Nishida nilpotence}

We now concentrate 
on the first classical theorem, that is Nishida nilpotence.
Its proof relies on the following
three key ingredients: the fact -- due to Serre -- that all elements 
in $\pi_*(S)$ are torsion for $* \neq 0$,
a careful study of power operations for
$H_{\infty}$-ring spectra and the Kahn-Priddy
theorem \cite{KP}. (We do not discuss the slightly different proof for those
elements annihilated by some prime $p$ already,  
which provides a sharper bound.) 
One might ask if there is a motivic generalization
of Nishida's theorem:

\begin{question}
For which fields $k$ and in which bidegrees
the following is true
for every $f \in \pi_{**}(S)$: If
$r \cdot f=0$ for some $r \in \Z$, then
$f$ is nilpotent?
\end{question}

If this was true for some field in all bidegrees, then every
non-nilpotent element $f \in \pi_{**}(S)$ would be
non-torsion, that is $r \cdot f \neq 0$ for all $r \in \Z$.
But for $k=\C$ both $\eta$ and $\mu_9$ are in the 
$\tau$-local region as described in \cite{GhI},
so they must be torsion by Levine's comparison theorem \cite{L}
and Serre's classical finiteness result.
Hence for $k=\C$ the answer is no in arbitrary
bidegrees, and yes when restricted to weight zero.
Next, one might look at $k = \R$ where we only have 
$(1 - \epsilon) \cdot \eta = 0$. This does
not contradict motivic Nishida nilpotence, except if
we generalize from $\Z$- to $GW(k)$-torsion.
(Note that there is an $\Z/2$-equivariant version of Nishida's
theorem, see \cite{Ir}. This might be a hint that everything that 
goes wrong with nilpotence goes wrong over $\C$ 
in some sense.)
Also, the elements $\eta$ and $\mu_9$ are $2$-torsion,
so the conjecture might be true for elements
annihilated by (products of) powers of odd primes. 
Indeed, the recent computaions by Stahn \cite{St}
at odd primes for $k=\C, \R$ have not yet led to new non-nilpotent
self-maps.
\medskip

According to Guillou and Isaksen \cite{GuI}, 
over $\C$ completed at $2$, besides the units in $\pi_{0,0}(S)$ 
there is an infinite family of non-nilpotent
elements $\mu_{8k+1} \in \pi_{8k+1,4k+1}(S)$ 
detected by $P^kh_1$, and starting
with $\mu_1=\eta$ and $\mu_9$. They further claim
(compare also \cite{GhI}) that there are other
families of non-nilpotent elements, e.g. one
starting with an element $d_1 \in \pi_{32,18}(S)$. 
According to Isaksen, the element $d_1$
even lifts to an element over $\R$, which consequently
is non-nilpotent either, and we have $4d_1=0$ over $\R$.
The element $d_1$ lives in the motivic 
four-fold Toda bracket $<\eta,\sigma^2,\eta,\sigma^2>$,
which is non-empty by Corollary \ref{d1int}.
We avoid to include these
unpublished results in our statements
and proofs, except for
Corollary \ref{d1int} on $d_1$ below. 

We also note that higher powers of $\nu$ and $\sigma$ lie in the
``not understood''-region. However, Isaksen's computations
show that they are nilpotent over $\C$. Based on this and
his computations over $\R$, Isaksen conjectures
that $\nu^4=0$ and $\sigma^4=0$ over any base
scheme.  When restricting to simplicial degree zero,
we have to study nilpotence for $K^{MW}_*(k)$
by Morel's theorem \cite[Theorem 6.2.1]{Mo},
at least if $k$ is perfect and $char(k) \neq 2$.
In this case, the appendix of Marcus Zibrowius shows
that Nishida nilpotence holds in non-positive degrees
-- that is non-negative degrees in the indexing of $K^{MW}_*(k)$
-- and in all degrees if and only if $k$ is formally real.

\medskip

So what goes wrong when trying to translate
the classical proof for Nishida nilpotence
to the motivic case?
Constructing extended powers and power operations in
stable motivic homotopy theory is possible.
For this, it is convenient to
use motivic symmetric spectra as introduced by Hovey and Jardine.
In what follows, we will work with motivic strictly commutative ring spectra
which are in particular motivic $H_{\infty}$-spectra.
That is, for any motivic strictly commutative
ring spectrum ${\bf E}$ we have maps
$$\xi_j: D_j{\bf E}:=(E \Sigma_j)_+ \wedge _{\Sigma_j} {\bf E}^{\wedge j}
\to {\bf E}$$ for all $j \geq 1$ such that
the diagrams of \cite[Definition I.3.1]{M2}
commute in $\SH(k)$.
These maps are given by $\xi_j=\beta_{j,0}$
with \cite[Remark I.2.6]{M2} applied to the motivic
setting. Virtually everything else
in \cite[sections I.2+3]{M2}
and the definitions of power operations, that is
\cite[Definitions I.4.1 and I.4.2]{M2},
then easily extends to the motivic setting.
We are mainly interested in the case ${\bf E}=S$, of course.
 
Having constructed these motivic power operations, we wish to
proceed similarly to \cite{M1} and \cite[sections II.1+2]{M2}.
For this, we need to prove a motivic version of
the Kahn-Priddy theorem (see below). The original proof of \cite{KP} uses the 
Barratt-Priddy-Quillen Theorem $B\Sigma_{\infty}^+ \simeq QS^0$
and various computations for (co-)homology with finite coefficients.
For ``geometric'' classifying spaces $BG$,
Voevodsky \cite[section 6]{V} has computed $H^{**}(BG,\Z/l)$
for $G=\Sigma_t, \mu_t$, which could be useful here.
In an unpublished draft some time ago, Morel conjectured a motivic 
version of the Barratt-Priddy-Quillen theorem
and pointed out a possible relationship with Serre's 
splitting principle for \'etale algebras. 

\begin{definition}
We say that the unstable motivic Kahn-Priddy theorem holds
for a field $k$ at the prime $p$ if the
motivic version of the map $\tau_p=h_p$ of
\cite[Definition II.1.4]{M2}, \cite[Lemma 1.6]{M1},
see Definition \ref{deftau} below,  
induces an epimorphism
$$(\tau_p)_*:\pi_{**}(D_pS^{q,w}) \to \pi_{**}(S^{pq,pw})$$
in bidegrees $(r,*)$, provided $r$ lies in the classical
Kahn-Priddy range of \cite[Theorem II.2.8]{M2}.
\end{definition}

Similar to the classical case \cite{KP}, this is an unstable variant 
which is related to a stable one.
The stable variant would predict epimorphisms
(everything localized at $p$)
$$(\tau_p)_*:\pi_{**}(\Sigma_T^{\infty}B\Sigma_p) \to \pi_{**}(S)$$
in a certain range.
There are many variants of the Kahn-Priddy theorem,
both for the statement and for the proof.
See e.g. Segal, Caruso-Cohen-May-Taylor,
L\"offler-Ray...

\medskip
 
Let us explain the map $\tau_p$ and some of its properties
in the motivic setting. It is possible to prove
the following two propositions for arbitrary monoidal
model categories tensored over pointed simplicial sets.
In particular, they are true for ``global'' model structures,
that is before carrying out motivic localizations.

We start with a motivic variant of \cite[Theorem II.1.1]{M2}.
We fix positive integers $j$ and $k$ and motivic (symmetric) $T$-spectra
$Y_1,...,Y_k$. Let $Z:= Y_1 \vee ... \vee Y_k$, and let
$\nu_i:Y_i \to Z$ be the inclusions. For a partition
$J=(j_1,...,j_k)$ of $j=j_1+...+j_k$, let
$f_J$ denote the composite
$$D_{j_1}Y_1 \wedge ... \wedge D_{j_k}Y_k \stackrel{D_J(\nu_J)}{\to}
D_{j_1}Z \wedge ... \wedge D_{j_k}Z \stackrel{\alpha_J}{\to}D_jZ.$$
Here $D_J(\nu_J):=D_{j_1}\nu_1 \wedge ... \wedge D_{j_k}\nu_k$,
and $\alpha_J$ is induced by the multiplication of the commutative
ring spectrum, and is well-defined by the motivic analogue of 
\cite[Lemma I.2.8]{M2}. 

\begin{proposition}\label{thm1.1}
In the above situation, the wedge sum
$$f_j: \bigvee_{J} D_{j_1}Y_1 \wedge ... \wedge D_{j_k} Y_k \to D_j Z$$
of the maps $f_J$ is a stable motivic equivalence.
\end{proposition}
\proof
This is similar to \cite[Theorem II.1.1]{M2}, using that motivic
spectra are tensored over simplicial sets. The map
$i \wedge 1$ corresponding to the one considered at the 
end of the proof of loc. cit. is a weak equivalence
already for global model structures, that
is before motivic localization.
\qed

For a fixed partition $J=(j_1,...,j_k)$ of $j$, let 
$g_J:D_j(Z) \to D_{j_1}Y_1 \wedge ... \wedge D_{j_k} Y_k$
be the $J$th component of $f_j^{-1}$. We now restrict to
the special case $Y:=Y_1=...=Y_k$, hence $Z=\bigvee_{i=1}^k Y$.
Let $\Delta: Y \to \bigvee_{i=1}^k Y = \prod_{i=1}^k Y$ be the diagonal map.
 
\begin{definition}\label{deftau}
For $J$ and $Y$ as above, we let $\tau_J$ be the composition
$$D_j Y \stackrel{D_j \Delta}{\to} D_j  \bigvee_{i=1}^k Y 
\stackrel{g_J}{\to}  D_{j_1}Y \wedge ... \wedge D_{j_k} Y.$$
If $J=(1,...,1)$ and hence $k=j$, we set 
$$\tau_j:=\tau_J: D_j Y \to D_1 Y \wedge ... D_1 Y = Y^{\wedge j}.$$
\end{definition}

The following is a motivic generalization of 
\cite[Lemma 1.6]{M1}, \cite[Corollary II.1.8]{M2}:

\begin{proposition}\label{coro1.8}
If $r=p^iv$ with $p$ prime, $i \geq 1$, and $v$ prime to $p$,
then $D_p(r):D_p Y \to D_p Y$ can be written as 
$p^i \lambda + (p,r-p)\iota_p \tau_p$ 
for some self-map
$\lambda$, where $\iota_p:Y^{\wedge p} \to D_p(Y)$ is
the canonical inclusion, where $(p,r-p):=\binom{r}{p}$. 
\end{proposition}
\proof
Similar to \cite[Corollary II.1.8]{M2}, using
Proposition \ref{thm1.1} and other results above.
\qed

From this, we can deduce the following key
result, which is is essentially a motivic generalization 
of \cite[Corollary II.2.4]{M2} and \cite[Theorem 3.8]{M1}.
(More precisely, setting $n=1$ in the latter
these two results are equivalent except that
\cite{M2} considers $\alpha \in E_{r}(D_pS^q)$
and \cite{M1} considers $\tilde{y} \in \pi_{pr+t}(D_pS^r)$.
However, we later restrict to $E=S$ anyway.)
The proof of \cite{M2} uses additivity formulae 
for power operations, which is some refinement
of the classical formulae for powers of sums.
We follow the proof hinted at in \cite{M1} instead,
see also the remark after \cite[Corollary II.2.5]{M2}.

\begin{proposition}\label{coro2.4}
Fix a strictly commutative (or $H_{\infty}$) motivic
$p$-local ring spectrum ${\bf E}$ and 
$x \in \pi_{q,w}{\bf E}={\bf E}_{q,w}$.
If $p^i \cdot x=0$ for some positive integer $i$, then
$p^{i-1} \cdot (\tau_p)_*(\alpha) \cdot x^{p+1}=0$ for all 
$\alpha \in \pi_{p(q,w)+(s,t)}(D_pS^{q,w})$.
\end{proposition}
\proof
As some details in the proof of 
the classical analogue \cite[Theorem 3.8]{M1}
are omitted, we will include these
here. Similar to loc. cit, we really show that 
$p^{i-1}(pz + x^py)=0$ for some $z:S^{p(q,w) +(s,t)}\to {\bf E}$
with $y:= (\tau_p)_*(\alpha)$.
(Following the tradition of \cite{M1}, our notation 
here and below omits $\Sigma^{\infty}_T$, and ocassionally does not
distinguish between a stable map and its unstable representative.) 
Multiplying this equality with $x$, 
we deduce that $p^{i-1} x^{p+1} (\tau_p)_*(\alpha)=0$
as claimed.
To show the equality above, we need to compute
$x^p \cdot y$ which is given by the following
composition:
$$S^{p(q,w)+(s,t)} \stackrel{\alpha}{\to} D_p(S^{q,w}) \stackrel{\tau_p}
{\to}S^{p(q,w)} 
\stackrel{\iota_p}{\to}D_p(S^{q,w}) \stackrel{D_p(x)}{\to} D_p({\bf E})
\stackrel{\xi_p}{\to} {\bf E}.$$
Now we multiply this composition with
$p^{i-1}$ and apply Proposition \ref{coro1.8} with
$r=p^i$ (hence $\binom{r}{p}$ is divisible by
$p^{i-1}$)  to  $p^{i-1} \cdot \iota_p \circ \tau_p$.
Consequently, the latter equals (possibly up to a unit) $D_p(p^i)-p^ig$ 
for some $g$. Hence $p^{i-1} \cdot \iota_p \circ \tau_p$
can be written as the difference of two maps, one being
already of the requested form and the other one given by
$$ S^{p(q,w)+(s,t)} \stackrel{\alpha}{\to} D_p(S^{q,w}) \stackrel{D_p(p^i)}{\to}
D_p(S^{q,w}) \stackrel{D_p(x)}{\to} D_p({\bf E})
\stackrel{\xi_p}{\to} {\bf E}.$$ 
But that one is zero as $D_p(p^i)D_p(x)=D_p(p^ix)=0$
by assumption.
\qed

\begin{theorem}\label{main}
Assume that the unstable motivic Kahn-Priddy theorem holds
for a field $k$ at a prime $p$. Then for $s \geq 1$,
any element $x \in \pi_{s,t}(S)$
over $k$ annihilated by a power of $p$ is nilpotent.
\end{theorem}
\proof
We may argue as in \cite[Theorem II.2.9]{M2}.
Namely, using the unstable motivic Kahn-Priddy theorem
we see that $x=(\tau_p)_*(\alpha)$ for some
$\alpha \in \pi_{p(q,w)+(s,t)}(D_pS^{q,w})$
with $q$ large enough and arbitrary $w$. 
More precisely, as in \cite[Theorem 2.9]{M2} we 
have a factorization $q=m \cdot s$ for a suitable integer
$m$, and then we may take $w=m \cdot t$. 
Now thanks to our assumption $p^ix=0$,
we can apply Proposition \ref{coro2.4}
with ${\bf E}=S$ to $\alpha$ and $x^m$ (rather 
than $x$). This yields the desired
equality $p^{i-1}x^{1+m(p+1)}=0$.
We then may conclude via descending induction 
on $i$. (Note that for odd primes $p$,
if $q-w$ is odd and $w$ is even or $k=\C$,
then we already have $x^2=0$ by the graded
commutativity of \cite[Proposition 2.5]{DI3}.)
\qed

As both $\eta$ and $\mu_9$ are $2$-torsion over
$\C$ (more generally $\eta$ is if $k$ is not formally 
real, as then $1 - \epsilon = 2$) 
and both are not nilpotent, the theorem above puts
some restrictions on possible motivic
Kahn-Priddy theorems.
Let us look at what goes wrong for $\eta \in \pi_{1,1}(S)$.

\begin{theorem}\label{mainappl}
The unstable motivic Kahn-Priddy theorem does not
hold at $p=2$. In particular:
\begin{itemize}
\item
The map
$$(\tau_2)_*:\pi_{9,2w+1}(D_2 S^{4,w}) \to \pi_{9,2w+1}(S^{8,2w})$$
for $w=1,2,3$ and $k$ not formally real
(e.g. $k=\C$) is not a $2$-local epimorphisms.
\item
The map
$$(\tau_2)_*:\pi_{1161,645}(D_2 S^{576,320}) \to \pi_{1161,645}(S^{1152,640})$$
for $k=\C$ is not a $2$-local epimorphism, provided
$\mu_9$ lifts to an unstable element in $\pi_{1161,645}(S^{1152,640})$.
\end{itemize}
\end{theorem}
\proof
By Proposition \ref{coro2.4} above, we know
the following: for any $\alpha \in \pi_{p(q,w)+(s,t)}(D_2S^{q,w})$
and any $x \in \pi_{q,w}(S)$ with $2 \cdot x =0$,
we have $(\tau_2)_*(\alpha) \cdot x^3=0$.
Now we need to find some $\alpha$ with
$(\tau_2)_*(\alpha)=\eta$. 
In the classical case (corresponding to $w=0=t$),
the unstable Kahn-Priddy theorem tells us that
there is some $\alpha \in \pi_{2 \cdot 4 +1}(D_2S^4)$
which is mapped to $\eta_{top} \in \pi_9(S^8)$
under $(\tau_2)_*$. Hence we may choose $q=m\cdot 1=4$
and apply Proposition \ref{coro2.4} to $x=\eta_{top}^4$.
The corresponding motivic map $(\tau_2)_*$
must not have a preimage $\alpha$ in certain
bidegrees $(9,w)$
for $\eta$. Namely, the claim about the first map
follows from Theorem \ref{main}, the 
above properties of $\eta$ and the above discussion
on $(\tau_2)_*$. For the second claim, the fact
that $\eta$ lifts to $\pi_{3,2}(S^{2,1})$,
$2$ lifts to $\pi_{1,0}(S^{1,0})$
and $\sigma$ lifts to $\pi_{15,8}(S^{8,4})$ \cite[section 4]{DI3}
together with some explicit construction
of Toda brackets makes it plausible
that $\mu_9$ lifts to an element in $\pi_{20+s,10+w}(S^{11+s,5+w})$
for small $s$ and $w$ similar to the classical case,
although we still lack a proof of this.
(The classical proof
of \cite[Theorem II.2.9]{M2} seems to use
Freudenthal's suspension theorem without mentioning it.
So far we have a motivic version of this
only with respect to the simplicial sphere, see
\cite[Theorem 6.61]{Mo2}.) If this holds,
the $\mu_9$ is also an element in $\pi_{1161,645}(S^{1152,640})$.
Now an argument similar to the one above for $\eta$ 
leads to the claim about the second map, choosing 
$(q,w)= 64 \cdot (9,5)=(576,320)$.
\qed
 
We leave it to the reader to deduce similar statements
using other non-nilpotent elements mentioned earlier.
Note that the element $d_1 \in \pi_{32,18}(S)$,
see also Lemma \ref{d1int} below,
shows that the motivic Kahn-Priddy theorem also fails
over $\R$, but we do not know in which precise bidegrees
as we do not dispose of a specific unstable lift of $d_1$ yet.
Still, it seems likely that the motivic
Kahn-Priddy theorem holds in weight $w=0$.

\medskip

Over $\R$, there are further non-nilpotent elements 
in $\pi_{**}(S)$ which do not exist over $\C$, of course:

\begin{lemma}
For $k=\R$, the elements $\epsilon \in \pi_{0,0}(S)$
and $\rho=\rho_{-1}=[-1] \in \pi_{-1,-1}(S)$ are not nilpotent.
\end{lemma}
\proof
We have $\epsilon^2=1$. Concerning $[-1]$, one knows that
there is a graded ring epimorphism $K_*^M(\R) \to \Z/2[t]$
given by $[-1] \mapsto t$. Hence $[-1]$ is not nilpotent
in $K_*^M(\R)$, and consequently not in $K_*^{MW}(\R)$.  
\qed
  
Note that $[-1]$ is detected by ${\bf MGL}_{**}$ because 
${\bf MGL}_{-1,-1}\cong K_1(\R)$ \cite[Theorem 6.4.5]{Mo}. 
In this sense, it behaves more classical
\cite{DHS} than $\eta$ which is not detected.

Finally, for $k=\C$ we have an element $\tau \in
\pi_{0,1}(L_{S/2}S)$, sometimes also denoted by $\theta$,
where $L_{S/2}S$ denotes motivic Bousfield localization with
respect to the mod-$2$ Moore spectrum $S/2$.
The element $\tau$ can be explicitely constructed using inverse systems
of roots of unity, see \cite[p.81/82]{HKO1}.
It is not known if this element 
lifts to an element in ``integral'' $\pi_{0,-1}(S)$.
This is a special case of the following more
general problem:

For general base fields, the abutment of the motivic ASS or ANSS 
might be quite different from $\pi_{**}(S)$.
Recall that for Morel's \cite{Mo} computation in simplicial
degree $0$ these spectral sequences are not necessary,
and that the article \cite{OO} on computations
in simplicial degree $1$ contains techniques
to get rid of completions in some cases, see also Lemma
\ref{mu9int} below.
Even for $k=\C$, the situation is more complicated
than in classical stable homotopy theory. 
Over $\C$, the motivic ASS converges to 
the nilpotent completion $\pi_{**}(S^{\wedge}_{\bf H})$
\cite[Corollary 6.15]{DI2}, \cite{HKO1}, where
${\bf H}$ denotes the  motivic Eilenberg-MacLane spectrum
for $\Z/p$ with $p$ a fixed prime. Furthermore, it is shown in 
\cite[Theorem 6]{HKO1} that the nilpotent completion 
$S^{\wedge}_{\bf H}$ has the same bigraded homotopy groups
as the Bousfield localization $S^{\wedge}_p:=L_{S/p}S \simeq holim \, S/p^n$
\cite[Lemma 18]{HKO3}.
Assuming this, it still remains to understand the map
$$\pi_{**}(S) \to \pi_{**}(S^{\wedge}_p).$$
For this, one may use a motivic variant \cite{HKO1} of
the usual short exact sequence of Bousfield
\cite[Proposition 2.5]{Bo}. This tells us
that $\pi_{s,w}(S^{\wedge}_p) \cong \pi_{s,w}(S)^{\wedge}_p$
if $\pi_{s,w}(S)$ and  $\pi_{s-1,w}(S)$ are finitely
generated abelian groups, but in general the situation
is more complicated. For the motivic
ANSS, similar problems do occur. The most convenient approach 
is to apply a motivic variant of the chromatic square,
as in the following example.

\begin{lemma}\label{mu9int}
For $k=\C$, the canonical map
$$\pi_{9,5}(S) \to \pi_{9,5}(L_{S/2}S)$$
is an isomorphism. In particular, the element
$\mu_9$ of \cite{GuI}, \cite{AM} has a unique integral lift.
\end{lemma}
\proof
We consider the restriction of the motivic arithmetic
square of \cite[Appendix A]{OO} at the prime $2$,
that is for the motivic localization functors $L_{S_{\Q}}$
and $L_{S/2}$. As we know by computations of
\cite{DI2} and \cite{HKO1} that $\pi_{9,5}(L_{S/2}S)$
and $\pi_{8,5}(L_{S/2}S)$ are finite $2$-torsion,
the associated long exact sequence degenerates
to the isomorphism 
$$\pi_{9,5}(S) \stackrel{\simeq}{\to} \pi_{9,5}(L_{S/2}S) \oplus 
\pi_{9,5}(L_{S_{\Q}}S).$$
By Morel's theorem \cite{Mo}, see also \cite{CD}, together with
the fact that $H^{-9,-5}_{mot}(Spec(\C))=0$, we know that 
$\pi_{9,5}(L_{S_{\Q}}S)=0$.  
\qed

Observe that similar arguments apply to many other elements
over $k=\C$, e.g. also to the higher $\mu_{8k+1}$ and to $d_1$, 
but not to $\tau$. 
Using the the recent work of Ananyevskiy-Levine-Panin \cite{ALP}, 
we also get results over other base fields:

\begin{corollary}\label{d1int}
For $k=\R$, the element $d_1 \in \pi_{32,18}(L_{S/2}S)$
lifts to a unique non-nilpotent element in $\pi_{32,18}(S)$.
\end{corollary}
\proof
This is similar to the previous Lemma. One has to replace 
the computations of \cite{DI2} by the recent unpublished ones of 
Isaksen in bidegrees $(32,18)$ and $(31,18)$,
and Morel's theorem by its refinement
provided in \cite{ALP}
which implies the rational vanishing in the required degrees.
\qed

Hence, also for real fields there are other
examples than $\eta$ for non-nilpotent elements!
Finally, note that if we can
not establish an integral lift for an element
in the completion, then Theorem \ref{mainappl}
still holds with localization replaced by completion.

\section{Relating the motivic sphere spectrum to the Witt spectrum}

We now discuss an interesting relationship between
the $\eta$-local motivic sphere spectrum
and Balmer's $4$-periodic Witt groups, represented by ${\bf KT}$
\cite[Theorem 5.8]{Ho},\cite{ST}.
For this, consider the unit map $S \to {\bf KO}$
for hermitian $K$-theory, 
sometimes also denoted by ${\bf BO}$ or
${\bf KQ}$ rather than ${\bf KO}$.
Together with the equivalence
${\bf KO}[\eta^{-1}] \simeq {\bf KT}$,
it induces a map ${\bf KO}[\eta^{-1}] \to {\bf KT}$,
which is the unit map for the naive ring spectrum
${\bf KT}$. Observe that the equivalence 
${\bf KO}[\eta^{-1}] \simeq {\bf KT}$ can be deduced
from \cite[Theorem 4.4]{RO} arguing as in \cite[sections 4 and 5]{Ho}.
The following question
assumes that $4 \cdot (1 - \epsilon)^2 \cdot \sigma=0$
for the base field $k$. As explained above,
this fails for subfields $k \subset \R$.
It seems reasonable to expect it is true if $-1$ is a sum of squares.  

\begin{question}\label{beautiful}
For any field $k$ with $char(k)\neq 2$
and $4 \cdot (1 - \epsilon)^2 \cdot \sigma=0$,
does the above map
$$u:S[\eta^{-1}] \to
{\bf KT}$$
in $\SH(k)$ factor through a map
$$u:S[\eta^{-1},\mu_9^{-1}] \to
{\bf KT}$$
and if so, what can we say about this map?
\end{question}

One might think of this question as an integral refinement
of the recent rational result of \cite{ALP},
which provides an answer
with rational coefficients. Morel's \cite{Mo}
computation of $\pi_{**}(S)$ in simplicial degree $0$
is also related to this question, of course.

\medskip

We have a complete answer to the question over the complex numbers.

\begin{theorem}\label{main2}
If $k=\C$, the unit map $u$ of Question \ref{beautiful}
factors through a map
$$u:S[\eta^{-1},\mu_9^{-1}] \to
{\bf KT}.$$
The latter map induces an
epimorphism $$u_{**}:\pi_{**}(S[\eta^{-1},\mu_9^{-1}]) \to 
\pi_{**}{\bf KT}$$
with the kernel being another copy of $\pi_{**}{\bf KT}$
shifted by bidegree $(7,4)$. 
\end{theorem}
\proof 
We show that $u_{**}(\mu_9)$ is invertible in $\pi_{**}({\bf KT})$.
For this, we first recall that the groups
${\bf KT}^{*,*}$ are $(4,0)$- and $(1,1)$-periodic \cite{Ho}, and 
we have $W(\C)=W^{0,0}(\C)$. We then consider the following commutative diagram

$\xymatrix{
\pi_{9,5}(S)  \ar[d]^{R_{\C}} \ar[r]^{u_{**}} & KO_{9,5}  
\ar[d]^{R_{\C}} \ar[r]^{{\simeq}\ \ \ \ \ \ \ \ \ } & KT_{9,5} \simeq W({\C}) \simeq {\Z}/2
\\
\pi_9 \ar[r]_{u_*^{top}\ \ \ \ \ \ } & KO^{top}_9 \simeq {\Z}/2
}$

where $u$ and $u^{top}$ are unit maps of ring spectra
and ${R_{\C}}$ denotes complex topological realization.
(There still seems to be no published proof
of the folklore theorem that $R_{\C}({\bf KO})={\bf KO}^{top}$.
The best way to prove this is probably to apply
the geometric description of ${\bf KO}$ established
in \cite{ST}.)
The element $\mu_9$ is mapped to the topological
$\mu_9$ in the topological stable stem $\pi_9$
and then to the non-zero element under
$u^{top}$ by \cite[Theorem 1.2]{Ad}.
As the square commutes, it is also mapped
to the non-zero element in $KT_{9,5}$, which is invertible
by the periodicity isomorphisms above.
Hence we obtain the desired factorization.
We now prove that $$ u_{**}: \pi_{**}(S[\eta^{-1},\mu_9^{-1}])
\to {\bf KT}^{-*,-*} $$ is an epimorphism.
By the Theorem of \cite{AM} above, 
the groups $\pi_{**}(S[\eta^{-1}\mu_9^{-1}])$
are $(4,0)$ and $(1,1)$-periodic, and so are the groups
${\bf KT}^{*,*}$. As $u_{**}$ is induced by the unit,
it maps the generator of $\pi_{0,0}
(S[\eta^{-1}\mu_9^{-1}])$ to the generator
of $W^{0,0}(\C)$. Moreover, we already saw that $u_{**}$ maps the invertible elements
$\eta$ and $\mu_9$ to invertible elements.
The claim now follows from the explicit structure
of the two bigraded rings.
\qed

\begin{remark}
In a preliminary version of this article, it was incorrectly
stated that the map on $\pi_{**}$ is an isomorphism
when replacing $S[\eta^{-1},\mu_9^{-1}]$ by
$S/(\sigma)[\eta^{-1},\mu_9^{-1}]$. 
This is obviously wrong, as Bachmann pointed out.
Note that in general the behaviour of motivic
homotopy groups when coning out elements is at least
as interesting as in the classical case, see e.g.
\cite[Proposition 9.3.2]{J}
\end{remark}

What can we say about the question for other base fields?
Assume that we have $4 \cdot (1 - \epsilon)^2 \cdot \sigma=0$,
or $2^t \cdot (1 - \epsilon)^{4-t} \cdot \sigma=0$ 
for some other $t$. 
Then the element $\mu_9$ exists, and it is non-nilpotent
-- in particular non-zero -- for any subfield $k \subset \C$
as the base change $\SH(k) \to \SH(\C)$ preserves
multiplication, Toda brackets and motivic Hopf elements.
Hence the factorization in the question
follows if we can show that
$\mu_9$ is mapped to a unit (this might follow
for subfields $k \subset \C$ as above). Now
is $u_*: \pi_{**}(S[\eta^{-1},\mu_9^{-1}])
\to {\bf KT}^{-*,-*}$
an epimorphism in all bidegrees? By periodicity,
this reduces to show the following. First, $u_*$ is an isomorphism in
bidegree $(0,0)$. This corresponds to Morel's theorem \cite{Mo}
before inverting $\mu_9$.
Second, the groups $\pi_{s,0}(S[\eta^{-1},\mu_9^{-1}])$
have to be studied for $s=1,2,3$. This will be the hard part, of course.
The computations of $\pi_{1,*}(S)$ of \cite{OO} for 
``low dimensional fields'' fields and very recently by
R\"ondigs-Spitzweck-Ostvaer \cite{RSO2} 
and R\"ondigs \cite{R} for general base fields
might be useful here, as well as their study
of the behaviour of the unit map $S \to {\bf KT}$
on slices.  

Observe that looking at the real realization functor
does not help us, as it maps $\mu_{9}$ to zero.
We have not used the main theorem of \cite{ALP} 
when we proved the conjecture for $k = \C$. However, looking
at other fields it tells us at least that the odd
torsion before inverting $\mu_9$ and killing $\sigma$
is not too large. As the conjecture does not apply
to $k= \R$, the computations
of \cite{DI4} do not apply here, either.
The recent preprint \cite{GI2} contains computations 
for $\pi_{*,*}(S)[\eta^{-1}]$ over $\R$,
and even more recently we have
the theorem of Bachmann \cite{B} at odd primes.
Also, there is work in progress by R\"ondigs 
(now available, see \cite{R}) on
a related cell structure of ${\bf KT}$
over $\C$ which also relies on \cite{AM}
and implies that $\pi_{**}(S[\eta^{-1}])$
vanishes in simplicial degree $0$ and $1$. 
Finally, it remains to state
a precise conjecture relating $S$ and ${\bf KT}$
(or some other spectrum?)
for those fields not covered by Question \ref{beautiful}.
Part (2) of \cite[Theorem 1.1]{DI4} gives a
first hint on what kind of phenomena
might occur.

Jens Hornbostel, 
Bergische Universit\"at Wuppertal,
Fachgruppe Mathematik und Informatik,
Gau{\ss}strasse 20, 42119 Wuppertal,
hornbostel@math.uni-wuppertal.de.



\newpage

\section*{Appendix: Nishida Nilpotence in Milnor-Witt K-Theory}

\setcounter{theorem}{0}
\renewcommand{\thesection}{A}

{\em By Marcus Zibrowius}

\medskip






As in the main part of this paper, we say that Nishida nilpotence holds in a certain degree of a graded ring if all torsion elements in that degree are nilpotent.  Question~2.1 above asks over which fields and in which bidegrees Nishida nilpotence holds for \(\pi_{**}(S)\), where \(S\) is the motivic sphere spectrum.
In this short appendix, we answer this question for the zero-line \(\bigoplus_d \pi_{d,d}(S)\) as follows:

\begin{quote}
Let \(S\) denote the motivic sphere spectrum over a perfect field \(F\) of characteristic not two.  Nishida nilpotence holds in \(\pi_{d,d}(S)\) for all non-positive \(d\).  It holds for all \(d\) if and only if \(F\) is formally real.
\end{quote}

Our answer rests on Morel's concrete description of the zero-line as the Milnor-Witt K-ring of the base field:  \(\pi_{d,d}(S) \cong \KMW_{-d}(F)\) for any field \(F\) as above \cite[Thm~6.2.1]{Mo:NATO}.  Indeed, given this identification, the answer is immediately implied by the following result:

\begin{proposition}\label{app-prop}
Let \(F\) be a field of characteristic not two.  
Nishida nilpotence holds in all non-negative degrees of the Milnor-Witt K-ring \(\KMW_*(F)\).
It holds in all degrees of \(\KMW_*(F)\) if and only if \(F\) is formally real.
\end{proposition}

The remainder of this appendix constitutes a proof of this proposition.  

Throughout, \(F\) will denote a field of characterstic not two.  Recall that \(F\) is either formally real (i.\,e.\ there exists at least one ordering on \(F\)) or non-real (i.\,e.\ \(-1\) is a sum of squares in \(F\)), and that these two possibilities are mutually exclusive \cite[Thm~VIII.1.10]{Lam}.
For non-real fields, Nishida nilpotence fails because the non-nilpotent element \(\eta\in \KMW_{-1}(F)\) is torsion in this case.  Before we describe the situation in non-positive degrees in more detail, let us recall the structure of the Milnor-Witt K-ring in these degrees:
\[
  \KMW_n(F) \cong \begin{cases}
    GW(F) & \text{ for } n =0 \\
    \phantom{G}W(F) & \text{ for } n < 0\\
    \end{cases}
\]
The ring structure on \(\KMW_{\leq 0}(F)\) is determined by the ring structure on \(GW(F)\) and the fact that multiplication with the element \(\eta\in \KMW_{-1}(F)\) corresponding to the unit \(\langle 1\rangle\in W(F)\) induces the canonical projection \(GW(F)\twoheadrightarrow W(F)\) from degree~\(0\) to degree \(-1\) and the identity in lower degrees.  In particular, any homogeneous element of \(\KMW_{<0}\) can be written as \(\phi\eta^n\) for some \(\phi\in W(F)\) and some \(\eta \geq 0\). 

\begin{proposition}\label{app-prop-}
  All torsion in \(\KMW_{\leq 0}(F)\) is \(2\)-primary torsion. In degree zero, the nilpotent elements are precisely the torsion elements.  In negative degrees, the homogeneous nilpotent elements are precisely the elements of the form \(\phi\eta^n\) with \(\phi\in W(F)\) a torsion element of even rank.
\end{proposition}
The assertions of Proposition~\ref{app-prop} concerning non-positive degrees easily follow from this proposition in view of facts $(W_3)$ and $(W_4)$ below.
Proposition~\ref{app-prop-} itself is just a reformulation of the remaining following well-known facts:
\begin{enumerate}
\item[($GW_1$)] \label{app-fact-GW-torsion}
  All torsion in \(GW(F)\) is two-primary.
\item[($GW_2$)] \label{app-fact-GW-nilpotence}
The nilpotent elements in \(GW(F)\) are precisely the torsion elements.
\item[($W_1$)] \label{app-fact-W-torsion}
  All torsion in \(W(F)\) is two-primary.
\item[($W_2$)] \label{app-fact-W-nilpotence}
  The nilpotent elements in \(W(F)\) are precisely the torsion elements of even rank.
\item[($W_3$)] \label{app-fact-formally-real}
  When \(F\) is formally real, all torsion elements in \(W(F)\) are of even rank. 
\item[($W_4$)] \label{app-fact-nonreal}
  When \(F\) is non-real, \(W(F)\) is torsion. In particular, the unit \(\langle 1 \rangle \in W(F)\) is a torsion element of odd rank.
\end{enumerate}
Facts ($W_3$) and ($W_4$) are immediate consequences of Pfister's description of the torsion subgroup of \(W(F)\) as the kernel of the total signature homomorphism \cite[Thm~VIII.3.2]{Lam}: observe that the rank homomorphism \(W(F)\to \Z/2\) factors through any signature. Fact $(W_1)$ is likewise stated in \cite[Thm~VIII.3.2]{Lam}.  For $(W_2)$, see  eq.~(8.16) at the end of section VIII.8 in \cite{Lam} and observe that nilpotent elements in \(W(F)\) necessarily have even rank.  The corresponding statements $(GW_1)$ and $(GW_2)$ easily follow by considering the following cartesian square of rings:
\[\xymatrix{
    {GW(F)} \ar[d]_{\text{rank}}\ar@{->>}[r] & {W(F)} \ar[d]_{\text{rank}}\\
    {\Z} \ar@{->>}[r] & {\Z/2}
}\]


It remains to prove Proposition~\ref{app-prop} in positive degrees. The idea is to use Morel's description of the Milnor-Witt K-ring as a fibre product of graded rings, generalizing the cartesian square above.  Let \(\KM_*(F)\) denote the Milnor K-ring of \(F\),  and let \(\kM_*(F)\) denote the Milnor K-ring modulo two. Write \(I(F)\subset W(F)\) for the fundamental ideal, consisting of all elements of even rank, and let  \(I_*(F)\) denote the graded ring given by \(W(F)\) in negative degrees and by \(I(F)^n\) in degrees \(n>0\).  By the positive answer to one of the famous questions in \cite{Milnor}, the graded Witt ring associated with the powers of the fundamental ideal is naturally isomorphic to \(\kM_*(F)\) \cite{OVV}.  In particular, we have a graded ring homomorphism \(I_*(F)\to \kM_*(F)\).  Morel shows in \cite[Thm~5.3]{Morel} that the Milnor-Witt K-ring fits into the following cartesian square of graded rings:
\[\xymatrix{
    {\KMW_*(F)} \ar[d]\ar[r] & {I_*(F)} \ar[d]\\
    {\KM_*(F)} \ar[r] & {\kM_*(F)}
}\]
We briefly dwell on the lower left corner.

\begin{lemma}\label{app-lem1}
Every element \(\alpha\) of positive degree in \(\KM_*(F)\) has a power of the form \(\alpha^m = 
\{-1\}\gamma\), for some \(m>0\) and some \(\gamma \in \KM_*(F)\).
\end{lemma}
(Both \(m\) and \(\gamma\) depend on \(\alpha\), of course.  The braces \(\{-\}\) indicate the canonical isomorphism from the units of \(F\) to \(\KM_1(F)\), translating multiplicative into additive notation.)
\proof
Suppose first that \(\alpha\) is a generator of \(\KM_n(F)\) (\(n>0\)). Then \cite[Lemma~1.2]{Milnor} implies
\(
     \alpha^2 = \pm \{-1\}^{ n}  \alpha
\).
In general, we can write \(\alpha\) as 
\(
    \alpha = \alpha_1 + \dots + \alpha_k
\)
for certain generators \(\alpha_i\).  
Then \(\alpha^{k+1}\) is a sum of products of the \(\alpha_i\)s, and in each summand at least one of the \(\alpha_i\)s appears at least twice. Thus, 
\(
\alpha^{k+1} = \{-1\}^n \gamma
\)
for some \(\gamma\).
\qed

\begin{lemma}\label{app-lem2}
Every element of positive degree in \(\KM_*(F)\) that becomes nil\-potent in \(\kM_*(F)\) is already nilpotent in \(\KM_*(F)\).
\end{lemma}

\proof
Let \(\alpha\) be such an element of positive degree. 
By assumption, \(\alpha\) has a power of the form 
\(
\alpha^k = 2 \beta
\)
for some \(\beta\).
By Lemma~\ref{app-lem1}, \(\alpha\) also has some power of the form
\(
\alpha^m = \{-1\} \gamma
\)
for some \(\gamma\). So
\(
\alpha^{m+k} = 2 \beta \{-1\} \gamma
\).
This vanishes as \( 2\{-1\} = \{(-1)^2\} = 0\) in \(\KM_1(F)\).
\qed

We now prove Proposition~\ref{app-prop} in degrees \(n>0\).  Using the above cartesian square, we can write any element of \(\KMW_n(F)\)  as a pair \((\alpha,\phi)\) with \(\alpha\in\KM_n(F)\) and \(\phi\in I^n(F)\) such that the image of \(\phi\) in \(\kM_n(F)\) agrees with the reduction of \(\alpha\) modulo two.  
Let \((\alpha,\phi)\) be such an element, and assume it is torsion.  
Then in particular, \(\phi\) is an even-rank torsion element of \(W(F)\), hence nilpotent in \(I_*(F)\) by ($W_2$). \textit{A fortiori,} its image in \(\kM_*(F)\) is nilpotent.  So by Lemma~\ref{app-lem2}, \(\alpha\) is nilpotent in \(\KM_*(F)\).  Hence \((\alpha,\phi)\) is nilpotent in \(\KMW_*(F)\), as claimed.



Marcus Zibrowius,
Heinrich-Heine-Universit\"at D\"ussel\-dorf,
Mathe\-mati\-sches In\-sti\-tut,
Universit\"ats\-stra{\ss}e 1, 40225 D\"ussel\-dorf,
marcus.zibrowius@uni-duesseldorf.de

\end{document}